
\magnification1200
\input amstex.tex
\documentstyle{amsppt}

\hsize=12.5cm
\vsize=18cm
\hoffset=1cm
\voffset=2cm

\footline={\hss{\vbox to 2cm{\vfil\hbox{\rm\folio}}}\hss}
\nopagenumbers
\def\DJ{\leavevmode\setbox0=\hbox{D}\kern0pt\rlap
{\kern.04em\raise.188\ht0\hbox{-}}D}

\def\txt#1{{\textstyle{#1}}}
\baselineskip=13pt
\def\hf{{\textstyle{1\over2}}}
\def\a{\alpha}\def\b{\beta}
\def\d{{\,\roman d}}
\def\e{\varepsilon}

\def\G{\Gamma}
\def\k{\kappa}

\def\={\;=\;}

\def\zt{\zeta(\hf+it)}

\def\D{\Delta}

\def\R{\Re{\roman e}\,} 
\def\z{\zeta}

\def\hf{{\textstyle{1\over2}}}
\def\txt#1{{\textstyle{#1}}}

\font\tenmsb=msbm10
\font\sevenmsb=msbm7
\font\fivemsb=msbm5
\newfam\msbfam
\textfont\msbfam=\tenmsb
\scriptfont\msbfam=\sevenmsb
\scriptscriptfont\msbfam=\fivemsb

\font\ff=cmr8
\def\txt#1{{\textstyle{#1}}}
\baselineskip=13pt

\font\teneufm=eufm10
\font\seveneufm=eufm7
\font\fiveeufm=eufm5
\newfam\eufmfam
\textfont\eufmfam=\teneufm
\scriptfont\eufmfam=\seveneufm
\scriptscriptfont\eufmfam=\fiveeufm
\def\mathfrak#1{{\fam\eufmfam\relax#1}}

\font\tenmsb=msbm10
\font\sevenmsb=msbm7
\font\fivemsb=msbm5
\newfam\msbfam
     \textfont\msbfam=\tenmsb
      \scriptfont\msbfam=\sevenmsb
      \scriptscriptfont\msbfam=\fivemsb

  \def\rightheadline{{\hfil{\ff
  Mean square for the zeta-function and the divisor problem }\hfil\tenrm\folio}}

  \def\leftheadline{{\tenrm\folio\hfil{\ff
   A. Ivi\'c }\hfil}}
  \def\emptyheadline{\hfil}
  \headline{\ifnum\pageno=1 \emptyheadline\else
  \ifodd\pageno \rightheadline \else \leftheadline\fi\fi}

\topmatter
\title
ON THE MEAN SQUARE OF THE ZETA-FUNCTION AND THE DIVISOR PROBLEM
\endtitle
\author   Aleksandar Ivi\'c  \endauthor
\address
Aleksandar Ivi\'c, Katedra Matematike RGF-a Universiteta u Beogradu,
\DJ u\v sina 7, 11000 Beograd, Serbia.
\endaddress
\email {\tt
ivic\@rgf.bg.ac.yu,  aivic\@matf.bg.ac.yu} \endemail
\dedicatory
\enddedicatory
\abstract
{Let $\D(x)$ denote the error term in the Dirichlet
divisor problem, and $E(T)$ the error term in the asymptotic
formula for the mean square of $|\zt|$. If
$E^*(t) = E(t) - 2\pi\D^*(t/2\pi)$ with $\D^*(x) =
 -\D(x)  + 2\D(2x) - \hf\D(4x)$, then we obtain
 the asymptotic formula
$$
\int_0^T (E^*(t))^2\d t \;=\; T^{4/3}P_3(\log T) + O_\e(T^{7/6+\e}),
$$
where $P_3$ is a polynomial of degree three in $\log T$ with
positive leading coefficient. The exponent 7/6 in the error term
is the limit of the method.
}
\endabstract
\endtopmatter

\document
\head
1. Introduction
\endhead

As usual, let the error term in the classical Dirichlet divisor problem be
$$
\D(x) \;=\; \sum_{n\le x}d(n) - x(\log x + 2\gamma - 1),
\leqno(1.1)
$$
and
$$
E(T) \;=\;\int_0^T|\zt|^2\d t - T\left(\log\bigl({T\over2\pi}\bigr) + 2\gamma - 1
\right),\leqno(1.2)
$$
where $d(n)$ is the number of divisors of
$n$, $\z(s)$ is the Riemann zeta-function, and $ \gamma = -\G'(1) = 0.577215\ldots\,$
is Euler's constant. In view of F.V. Atkinson's classical explicit formula
for $E(T)$ (see [1] and [2, Chapter 15]) it was known long ago that
there are analogies between $\D(x)$ and $E(T)$. However,
instead of the error-term function $\D(x)$ it
is more exact to work with the modified
function $\D^*(x)$ (see  M. Jutila [5], [6] and T. Meurman [9]), where
$$
\D^*(x) := -\D(x)  + 2\D(2x) - \hf\D(4x)
= \hf\sum_{n\le4x}(-1)^nd(n) - x(\log x + 2\gamma - 1),
\leqno(1.3)
$$
which is a better analogue of $E(T)$ than $\D(x)$.
M. Jutila (op. cit.) investigated both the
local and global behaviour of the difference
$$
E^*(t) \;:=\; E(t) - 2\pi\D^*\bigl({t\over2\pi}\bigr),
$$
and in particular in [6] he proved that
$$
\int_0^T(E^*(t))^2\d t \ll T^{4/3}\log^3T.\leqno(1.4)
$$

On the other hand, we have the asymptotic formulas
$$
\eqalign{&
\int_0^x\D^2(y)\d y \;=\; Ax^{3/2} + F(x),\cr&
A = {1\over6\pi^2}\sum_{n=1}^\infty d^2(n)n^{-3/2},\quad F(x)
\ll x\log^4x,\cr}\leqno(1.5)
$$
and
$$
\eqalign{&
\int_0^T E^2(t)\d t \;= \;BT^{3/2} + R(T),\cr&
B = {2\over3}(2\pi)^{-1/2}\sum_{n=1}^\infty d^2(n)n^{-3/2},\quad R(T)
\ll T\log^4T,\cr}\leqno(1.6)
$$
and an analogous formula to (1.5) holds if $\D(x)$ is replaced by $\D^*(x)$,
since Lemma 2 below is the analogue of the truncated Vorono{\"\i} formula
for $\D(x)$, and a full Vorono{\"\i} formula for $\D^*(x)$ follows from
the corresponding Vorono{\"\i} formula for $\D(x)$
(see e.g., [2, Chapter 3]) and (1.3). Note
that we have
$$\sum_{n=1}^\infty d^2(n)n^{-3/2} = {\z^4(3/2)\over\z(2)}.
$$
The bounds in (1.5) and (1.6) are at present the sharpest ones known, and the
history is given in [3, Chapter 2]. The bounds with $\log^5$ instead of $\log^4$
were independently obtained by Y. Motohashi and T. Meurman [10]. The actual
bound in (1.5) is due to E. Preissmann [11], and the bound in (1.6) was
noticed independently by Preissmann (loc. cit.) and the author [3].
This follows by using a variant of the so-called Hilbert's inequality
(see [3], eq. (2.101)]). Further results on $F(x)$ were obtained by Lau and
Tsang (see [7] and [13]). In particular, in [13] K.-M. Tsang has shown that,
for almost all $x$ and a suitable constant $\k$,
$$
F(x) \= -{1\over4\pi^2}x\log^2x + \kappa x\log x + O(x),
\leqno(1.7)
$$
and he conjectured that (1.7) holds unconditionally. Thus (1.5) seems to
be close to best possible, and very likely the same holds for (1.6).

\medskip
For an extensive discussion of $E(T)$
see the author's monographs [2], [3].
The significance of (1.4) is  that it shows that $E^*(t)$ is in
the mean square sense much smaller than either $E(t)$ or
$\D^*(t)$. It is expected that $E(t)$ and
$\D^*(t)$ are `close' to one another in order, but this has never
been satisfactorily established, although M. Jutila [5] obtained
significant results in this direction. It is also conjectured that
$E^*(T) \ll_\e T^{1/4+\e}, \D^*(x) \ll_\e x^{1/4+\e}$, the latter
being a consequence
of the classical conjecture  $\D(x)\ll_\e x^{1/4+\e}$ in the Dirichlet divisor
problem. In fact,Y.-K. Lau and K.-M. Tsang [8] recently
proved that $\a = \a^*$, where
$$
\a = \inf\Bigl\{a \;:\;\D(x) \ll x^a\,\Bigr\},\quad \a^* =
\inf\Bigl\{a^* \;:\;\D^*(x) \ll x^{a*}\,\Bigr\}.
$$

\medskip
In the first part of the author's work [4] the bound in (1.4) was complemented
with the new bound
$$
\int_0^T (E^*(t))^4\d t \;\ll_\e\; T^{16/9+\e};\leqno(1.8)
$$
neither (1.4) or (1.8) seem to imply each other.
Here and later $\e$ denotes positive constants which are arbitrarily
small, but are not necessarily the same ones at each occurrence,
while $a \ll_\e b$ means that the $\ll$--constant depends on $\e$.
In the second part of the same work (op. cit.) it was proved that
$$
\int_0^T |E^*(t)|^5\d t \;\ll_\e\; T^{2+\e},\leqno(1.9)
$$
and some further results on higher moments of $|E^*(t)|$ were obtained as well.
The aim of this note is to show that the bound in (1.4) is of the correct
order of magnitude. The result is the asymptotic formula, contained in

\bigskip
THEOREM 1. {\it We have
$$
\int_0^T (E^*(t))^2\d t \;=\; T^{4/3}P_3(\log T) + O_\e(T^{7/6+\e}),\leqno(1.10)
$$
where $P_3(y)$ is a polynomial of degree three in $y$ with
positive leading coefficient, and all the coefficients may be evaluated
explicitly.
}
\bigskip

From (1.10) and H\"older's inequality for integrals we obtain
the following

\medskip
{\bf Corollary 1}. {\it If $A\ge2$ is fixed, then}
$$
\int_0^T |E^*(t)|^A\d t \;\gg\; T^{1+A/6}(\log T)^{3A/2}.\leqno(1.11)
$$

\medskip
Note that there is a discrepancy between the bound in (1.11) (for $A=4$
and $A=5$), and the upper bounds in (1.8) and (1.9). The true order of the
integrals in (1.8) and (1.9) is difficult to determine. If, as usual,
$f(x) = \Omega(g(x))$ means that $\lim_{x\to\infty }f(x)/g(x) \ne 0$, then from Theorem 1
we immediately obtain

\medskip
{\bf Corollary 2}. {\it We have}
$$
E(T) = 2\pi\D^*\bigl({T\over2\pi}\bigr) + \Omega(T^{1/6}\log^{3/2}T).\leqno(1.12)
$$

\medskip
Thus equation (1.12) shows that $E(T)$ and $\D^*\bigl({T\over2\pi}\bigr)$ cannot be
too `close' to one another in order. A reasonable conjecture,
compatible with the classical conjectures
$E^*(T) \ll_\e T^{1/4+\e}, \D^*(x) \ll_\e x^{1/4+\e}$, is that one has
$$
E(T) = 2\pi\D^*\bigl({T\over2\pi}\bigr) + O_\e(T^{1/4+\e}).\leqno(1.13)
$$

\head
2. The necessary lemmas
\endhead

\bigskip
In this section we shall state three lemmas which are necessary
for the proof of Theorem 1. The first
two are Atkinson's classical explicit
formula for $E(T)$ (see [1] or e.g., [2] or [3])
and the Vorono{\"\i}-type formula for $\D^*(x)$, which is the analogue
of the classical truncated Vorono{\"\i} formula for $\D(x)$. The third
is an asymptotic formula involving $d^2(n)$.

\medskip
LEMMA 1. {\it Let $0 < A < A'$ be any two fixed constants
such that $AT < N < A'T$, and let $N' = N'(T) =
T/(2\pi) + N/2 - (N^2/4+ NT/(2\pi))^{1/2}$. Then }
$$
E(T) = \Sigma_1(T) + \Sigma_2(T) + O(\log^2T),\leqno(2.1)
$$
{\it where}
$$
\Sigma_1(T) = 2^{1/2}(T/(2\pi))^{1/4}\sum_{n\le N}(-1)^nd(n)n^{-3/4}
e(T,n)\cos(f(T,n)),\leqno(2.2)
$$
$$
\Sigma_2(T) = -2\sum_{n\le N'}d(n)n^{-1/2}(\log T/(2\pi n))^{-1}
\cos\left(T\log \Bigl( {T\over2\pi n}\Bigr) - T + {\txt{1\over4}}\pi\right),
\leqno(2.3)
$$
{\it with}
$$
\eqalign{\cr&
f(T,n) = 2T{\roman {arsinh}}\,\bigl(\sqrt{\pi n/(2T)}\,\bigr) + \sqrt{2\pi nT
+ \pi^2n^2} - {\txt{1\over4}}\pi\cr&
=  -\txt{1\over4}\pi + 2\sqrt{2\pi nT} +
\txt{1\over6}\sqrt{2\pi^3}n^{3/2}T^{-1/2} + a_5n^{5/2}T^{-3/2} +
a_7n^{7/2}T^{-5/2} + \ldots\,,\cr}\leqno(2.4)
$$
$$\eqalign{\cr
e(T,n) &= (1+\pi n/(2T))^{-1/4}{\Bigl\{(2T/\pi n)^{1/2}
{\roman {arsinh}}\,\Bigl(\sqrt{\pi n/(2T)}\,\Bigr)\Bigr\}}^{-1}\cr&
= 1 + O(n/T)\qquad(1 \le n < T),
\cr}\leqno(2.5)
$$
{\it and $\,{\roman{arsinh}}\,x = \log(x + \sqrt{1+x^2}\,).$}

\medskip
LEMMA 2 (see [2, Chapter 15]).{\it We have, for $1\ll N\ll x$},
$$
\D^*(x) = {1\over\pi\sqrt{2}}x^{1\over4}
\sum_{n\le N}(-1)^nd(n)n^{-{3\over4}}
\cos(4\pi\sqrt{nx} - {\txt{1\over4}}\pi) +
O_\e(x^{{1\over2}+\e}N^{-{1\over2}}).
\leqno(2.6)
$$
\medskip

LEMMA 3. {\it For $a > -\hf$ a constant we have}
$$
\sum_{n\le x}d^2(n)n^a = x^{a+1}P_3(\log x;a) + O_\e(x^{a+1/2+\e}),\leqno(2.7)
$$
{\it where $P_3(y;a)$ is a polynomial of degree three in $y$
whose coefficients depend on $a$,
and whose leading coefficient equals $1/(\pi^2(a+1))$. All the
coefficients of $P_3(y;a)$ may be explicitly evaluated}.

\medskip
{\bf Proof.} One can obtain (2.7) by partial summation from [2, eq. (14.30)],
but a direct proof seems in order. From the Perron inversion formula
(see e.g., the Appendix of [2]) we have, for $1 \ll T\ll x$,
$$
\sum_{n\le x}d^2(n) = {1\over2\pi i}\int_{1+\e-iT}^{1+\e+iT}
{\z^4(s)\over\z(2s)s}x^s\d s + O_\e(x^{1+\e}T^{-1}),\leqno(2.8)
$$
since $d^2(n)$ is generated by $\z^4(s)/\z(2s)$ (op. cit, eq. (1.106)).
The segment of integration in (2.8) is
replaced by segments joining the points $1+\e\pm iT, \hf+\e\pm iT$. There is a
pole of the integrand at $s=1$ of degree four, making a contribution
of $xP_3(\log x)$ by the
residue theorem, where $P_3(y)$ is a generic polynomial of degree
three in $y$. The leading coefficient of $P_3$ is found to be $1/ \pi^2$
(see [2, eq. (5.24)]). From $\zt \ll t^{1/6}$ the contribution of the vertical
segments $[\hf+\e \pm iT, 1+\e\pm iT]$
is found to be $\ll_\e x^{1+\e}T^{-1} + x^{1/2+\e}T^{-1/3}$.
By using the elementary bound
$$
\int_0^T|\z(\hf+\e+it)|^4\d t \ll_\e T
$$
for
the integral over the segment joining the points $\hf+\e\pm iT$, one arrives at
$$
\sum_{n\le x}d^2(n) = xP_3(\log x) + O_\e(x^{1/2+\e})\leqno(2.9)
$$
on taking $T = \sqrt x$ in (2.8),
from which (2.7) follows by partial summation. Note that we cannot rule out the
existence of zeros of $\z(s)$ in the region $\hf<\R s<1$. Therefore the error
terms in (2.9) and  (2.7) are essentially the best ones that we can obtain
without assuming some hypotheses on the zeros of $\z(s)$,
like e.g., the famous Riemann Hypothesis
(that  all complex zeros of $\z(s)$ are on $\R s = \hf$).
\head
3. The proof of  Theorem 1
\endhead
It is sufficient to prove the formula (1.10) with the integral
over $[T,\,2T]$ and then to replace $T$ by $T2^{-j}$ and to sum
all the results over $j = 1,2,\ldots\;$. We use Lemma 1 and Lemma 2
with $N=T$
to deduce that, for $T\le t\le2T$,
$$
E^*(t) = S_1(t) + S_2(t) +S_3(t), \leqno(3.1)
$$
where
$$
\eqalign{
&S_1(t) =\sqrt{2}\left({t\over2\pi}\right)^{1/4}
\sum_{n\le T}(-1)^nd(n)n^{-3/4}\Biggl\{e(t,n)\cos f(t,n)\cr&
\qquad\;- \cos (\sqrt{8\pi nt} - {\txt{1\over4}}\pi)\Biggr\}\d t\cr&
S_2(t) = -2\sum_{n\le c_1T}d(n)n^{-1/2}\left(\log {t\over2\pi n}\right)^{-1}
\cos\left(t\log {t\over2\pi n}-t+ {\pi\over 4}\right),\cr&
S_3(t) = O_\e(T^\e),\cr}\leqno(3.2)
$$
where $c_1 \,(>0)$ is to be found from the definition of
$N'$ in Lemma 1. We have
$$
\int_T^{2T}S^2_1(t)\d t \ll T^{4/3}\log^3T, \quad
\int_T^{2T}\{S^2_2(t)+S_3^2(t)\}\d t
\ll_\e T^{1+\e}.\leqno(3.3)
$$
The first bound is in fact embodied in the proof of (1.4), since it is $S_1(t)$
that makes the major contribution to the integral in question. This can be
seen if one follows the proof of (1.4) in [6] or in [2], and the same holds
for the estimate in (3.3) containing $S_2(t)$, while the estimate
containing $S_3(t)$ is trivial.
From (3.1)--(3.3) and the Cauchy-Schwarz inequality for integrals it
transpires that
$$
\int_T^{2T}(E^*(t))^2\d t = \int_T^{2T}S^2_1(t)\d t + O_\e(T^{7/6+\e}).
\leqno(3.4)
$$
By using the method of T. Meurman [10] the error term in (3.4) could be improved to
$O_\e(T^{1+\e})$, but it is the exponent 7/6 in (3.11) that cannot be improved
at present. Squaring out $S_1(t)$ it is seen that the integral on the right-hand
side of (3.4) equals
$$\eqalign{
&\sqrt{2\over\pi}\sum_{n\le  T}d^2(n)n^{-3/2}\int_T^{2T}t^{1/2}
\Biggl\{e(t,n)\cos f(t,n)
\cr&
- \cos (\sqrt{8\pi nt} - {\txt{1\over4}}\pi)\Biggr\}^2\d t\cr&
+ \sqrt{2\over\pi}\sum_{m\ne n\le T}(-1)^{m+n}d(m)d(n)(mn)^{-3/4}
\int_T^{2T}t^{1/2}\cdots,\cr}
\leqno(3.5)
$$
where the main contribution comes from the `diagonal' terms $m=n$, while
the terms with $m\ne n$ will contribute to the error term in (3.4).
The terms in the first sum in (3.5) for which $n>T^{2/3}$ trivially
contribute $\ll_\e T^{7/6+\e}$. The contribution of the second
sum in (3.5) (with $m\ne n$) is $\ll_\e T^{1+\e}$.
 This follows e.g., by the method of [2, Chapter 15]
or [3, Chapter 2] for the mean square of $E(t)$. Therefore (3.5) leads to
$$\eqalign{\int_T^{2T}(E^*(t))^2\d t &=
\sqrt{2\over\pi}\sum_{n\le T^{2/3}}d^2(n)n^{-3/2}\int_T^{2T}t^{1/2}
\Biggl\{e(t,n)\cos f(t,n)
\cr&
- \cos (\sqrt{8\pi nt} - {\txt{1\over4}}\pi)
\Biggr\}^2\d t + O_\e(T^{7/6+\e}),
\cr}
\leqno(3.6)
$$
and furthermore we use (see (2.5)) $e(t,n) = 1 + O(n/T)$ in (3.6) to see
that we may replace $e(t,n)$ by unity, with the ensuing error term absorbed by
the error term in (3.6). To manage the cosines in (3.6) we use the elementary
identity
$$
(\cos\a-\cos\b)^2 = 1 - \cos(\a-\b)+\{\hf\cos2\a +\hf\cos2\b - \cos(\a+\b)\}
$$
with $\a = f(t,n),\b = \sqrt{8\pi nt} - {\txt{1\over4}}\pi$.
By the first derivative test
(cf. [2, Lemma 2.1]) it is seen that the terms coming from curly braces
contribute $\ll T$ to (3.6). Likewise the values of $n$ in
$\cos(\a-\b)$ for which $T^{1/2} < n\le T^{2/3}$ contribute
$$
\ll  \sum_{n > T^{1/2}}d^2(n)n^{-3/2}T^{1/2}T^{3/2}n^{-3/2}
=  \sum_{n > T^{1/2}}d^2(n)n^{-3}T^2 \ll T\log^3T.
$$
For the terms $n\le T^{1/2}$ in $\cos(\a-\b)$ we use the series expansion
$$
\cos(\a-\b) = \cos(x_0+h) = \cos x_0  + \sum_{m=1}^\infty {h^m\over m!}
\cos(x_0 +\hf\pi m)\leqno(3.7)
$$
with ($a_1 = \txt{1\over6}\sqrt{2\pi^3}$, see (2.4))
$$
x_0 = a_1n^{3/2}t^{-1/2}, \; h = a_3n^{5/2}t^{-3/2} + a_5n^{7/2}t^{-5/2}+\ldots\,.
$$
We truncate the series over $m$ so that the tail, by trivial estimation,
 makes a total contribution which
is $\ll T$. This is possible because
$$
h \ll n^{5/2}T^{-3/2} \ll T^{-1/4}\qquad(T\le t \le 2T,\; 1 \le n \le T^{1/2}).
$$
The contribution of the finitely many remaining terms in the sum over $m$ in (3.7)
is estimated by the first derivative test, and it is
$$
\ll \sum_{n\le T^{1/2}}d^2(n)n^{-3/2}T^{1/2}n^{5/2}T^{-3/2}n^{-3/2}T^{3/2}
\ll T^{3/4}\log^3T.
$$
Finally it follows from (3.6) that we obtain, since
$1-\cos\gamma = 2\sin^2(\hf\gamma)$,
 $$\eqalign{
\int_T^{2T}(E^*(t))^2\d t &=
\sqrt{8\over\pi}\sum_{n\le T^{2/3}}d^2(n)n^{-3/2}\int_T^{2T}t^{1/2}
\sin^2(\hf a_1n^{3/2}t^{-1/2})\d t\cr&
+ O_\e(T^{7/6+\e}).\cr}\leqno(3.8)
$$
If we use the bound
$$
\sin^2(\hf a_1n^{3/2}t^{-1/2}) \;\ll\; n^3T^{-1}\qquad(T\le t \le 2T),
$$
then the contribution of the terms with $n < T^{4/15}$ in (3.8) is seen to be
$$\eqalign{&
\ll \sum_{n\le T^{4/15}}d^2(n)n^{-3/2}n^3T^{-1}\int_T^{2T}t^{1/2}\d t
\ll T^{1/2}\sum_{n\le T^{4/15}}d^2(n)n^{3/2} \cr&
\ll T^{1/2}(T^{4/15})^{5/2}\log^3T = T^{7/6}\log^3T.\cr}
$$
Hence (3.8) reduces to
$$\eqalign{
\int_T^{2T}(E^*(t))^2\d t &=
\sqrt{8\over\pi}\sum_{T^{4/15}\le n\le T^{2/3}}d^2(n)n^{-3/2}\int_T^{2T}t^{1/2}
\sin^2(\hf a_1n^{3/2}t^{-1/2})\d t\cr&
+ O_\e(T^{7/6+\e}).\cr}\leqno(3.9)
$$
In the integral above we make the change of variable
$$
 \hf a_1n^{3/2}t^{-1/2} = y,\; \d t = -\hf a_1^2n^3y^{-3}\d y.
$$
We set  $z = (4(y/a_1)^2T)^{1/3}$ so that, after changing the order of integration
and summation,  the main term on the right-hand side of (3.9) becomes
$$\eqalign{&
\sqrt{2\over\pi}\sum_{T^{4/15}\le n\le T^{2/3}}d^2(n)n^{-3/2}
\int_{{1\over2}a_1 n^{3/2}(2T)^{-1/2}}^{{1\over2}a_1n^{3/2}T^{-1/2}}
{a_1n^{3/2}\over y} \cdot{a_1^2n^3\over 2y^{3}}\cdot\sin^2y\d y\cr&
= \sqrt{1\over2\pi}a_1^3\int_{2^{-3/2}a_1 T^{-1/10}}^{{1\over2}a_1T^{1/2}}
\sum_{\max(T^{4/15},z)\le n\le \min(T^{2/3},2^{1/3}z)}d^2(n)n^3\cdot
{\sin^2y\over y^4}\d y.\cr}\leqno(3.10)
$$
The range of summation for $n$ is the interval $[z,\, 2^{1/3}z]$ if
$$
y\;\in\; J,\quad J := \left[{1\over2} a_1T^{-1/10},\,
{1\over2\sqrt{2}}a_1T^{1/2}\,\right].
$$
If we replace the interval of integration in the second integral in (3.10) by $J$,
then by using $|\sin y| \le\min(1,|y|)$ it follows that the error that is
made is $\ll_\e T^{7/6+\e}$.
Now we use Lemma 3 with $a=3$ to obtain that the  integral over $J$ equals,
with suitable constants $d_j$,
$$
\eqalign{&
\int_{{1\over2} a_1T^{-1/10}}^{{1\over2\sqrt{2}}a_1T^{1/2}}
\left(T^{4/3}y^{8/3}\sum_{j=0}^3d_j\log^j(y^{2/3}T^{1/3})
+ O_\e(T^{7/6+\e}y^{7/3})\right)\,{\sin^2y\over y^4}\d y\cr&
= T^{4/3}\int_{ {1\over2}a_1T^{-1/10}}^{{1\over2\sqrt{2}}a_1T^{1/2}}
\left(\sum_{j=0}^3d_j3^{-j}\log^j(y^2T)\right)\,{\sin^2y\over y^{4/3}}\d y
+ O_\e(T^{7/6+\e}).\cr}\leqno(3.11)
$$
Replacing the segment of integration in the integral on the right-hand
side of (3.11) by $(0,\,\infty)$, we make an error
which is $\ll_\e T^{7/6+\e}$. Namely, for
$0 < \a < 1 < \b,\; j = 0,1,\ldots\;$ fixed, we have
$$
\int_\a^\b\,{\sin^2y\over y^{4/3}}\log^j\d y =
\int_0^\infty\,{\sin^2y\over y^{4/3}}\log^j\d y
+ O(\a^{5/3}) + O(\b^{-1/3}\log^j\b),\leqno(3.12)
$$
where we used again $|\sin y| \le\min(1,|y|)$.
Hence the expression in (3.11) becomes, on using (3.12) with
$\a ={1\over2}a_1 T^{-1/10}$,
$\b = 2^{-3/2}a_1T^{1/2}$,
$$
T^{4/3}\sum_{j=0}^3b_j\log^jT +  O_\e(T^{7/6+\e})
$$
with some constants $b_j\,(b_3>0)$ which may be explicitly evaluated, and (1.10)
follows from (3.9) and (3.11). With some more effort one could improve the
error term in (1.10) to $O(T^{7/6}\log^CT)$, by using the method of
T. Meurman [10] (see also [3, Chapter 2]) and improving the error term
in Lemma 3 (see K. Ramachandra and A. Sankaranarayanan [12])
to $O(x^{a+1/2}\log^5x\log\log x)$, but the exponent 7/6 is  the limit
of the method.

\bigskip
\Refs
\bigskip

\item{[1]} F.V. Atkinson, The mean value of the Riemann zeta-function,
Acta Math. {\bf81}(1949), 353-376.

\item{[2]} A. Ivi\'c, The Riemann zeta-function, John Wiley \&
Sons, New York, 1985 (2nd ed. Dover, Mineola, New York, 2003).

\item{[3]} A. Ivi\'c, The mean values of the Riemann zeta-function,
LNs {\bf 82}, Tata Inst. of Fundamental Research, Bombay (distr. by
Springer Verlag, Berlin etc.), 1991.

\item{[4]} A. Ivi\'c, On the Riemann zeta-function and the divisor problem,
Central European J. Math. {\bf(2)(4)} (2004), 1-15, and II, ibid.
{\bf(3)(2)} (2005), 203-214.

\item{[5]} M. Jutila, Riemann's zeta-function and the divisor problem,
Arkiv Mat. {\bf21}(1983), 75-96 and II, ibid. {\bf31}(1993), 61-70.

\item{[6]} M. Jutila, On a formula of Atkinson, in ``{\it Coll. Math. Sci.
J\'anos Bolyai 34, Topics in classical Number Theory, Budapest 1981},
North-Holland, Amsterdam, 1984, pp. 807-823.

\item{[7]} Y.-K. Lau and K.-M. Tsang, Mean square of the remainder term
in the Dirichlet divisor problem, J. Th\'eorie des Nombres Bordeaux
{\bf7}(1995), 75-92.

\item{[8]} Y.-K. Lau and K.-M. Tsang, Omega result for the mean square of the
Riemann zeta-function, Manuscripta Math. {\bf117}(2005), 373-381.

\item{[9]} T. Meurman, A generalization of Atkinson's formula to $L$-functions,
Acta Arith. {\bf47}(1986), 351-370.

\item{[10]} T. Meurman, On the mean square of the Riemann zeta-function,
Quart. J. Math. (Oxford) {\bf(2)38}(1987), 337-343.

\item{[11]} E. Preissmann, Sur la moyenne quadratique du terme de reste
du probl\`eme du cercle, C.R. Acad. Sciences Paris S\'erie I {\bf306}(1988),
151-154.

\item{[12]} K. Ramachandra and A. Sankaranarayanan,
On an asymptotic formula of Srinivasa Ramanujan, Acta Arith. {\bf109}(2003), 349-357.

\item{[13]} K.-M. Tsang, Mean square of the remainder term in the Dirichlet
divisor problem II, Acta Arith. {\bf71}(1995), 279-299.

\endRefs

\enddocument

\bye